\documentclass[10pt]{amsart}
\usepackage{amsfonts}
\usepackage{amsmath}
\usepackage{amsthm}
\usepackage{amssymb}
\usepackage{latexsym}

\begin{document}
\newtheorem{cor}{Corollary}[section]
\newtheorem{lem}{Lemma}[section]
\newtheorem*{sublem}{Sublemma}
\newtheorem{prop}{Proposition}[section]
\newtheorem{ax}{Axiom}
\newtheorem{example}{Example}

\theoremstyle{definition}
\newtheorem{defn}{Definition}[section]
\theoremstyle{definition}
\newtheorem{thm}{Theorem}
\newtheorem*{conj}{Conjecture}
\newtheorem*{rem}{Remark}
\newtheorem{notation}{Notation}

\newenvironment{pf}{\proof}{\endproof}

\theoremstyle{remark}
\renewcommand{\thenotation}{} 

\numberwithin{equation}{section}
\setcounter{section}{-1}

\newcommand{\thmref}[1]{Theorem~\ref{#1}}
\newcommand{\secref}[1]{Sect.~\ref{#1}}
\newcommand{\lemref}[1]{Lemma~\ref{#1}}
\newcommand{\propref}[1]{Proposition~\ref{#1}}
\newcommand{\corref}[1]{Corollary~\ref{#1}}
\newcommand{\remref}[1]{Remark~\ref{#1}}
\newcommand{\er}[1]{(\ref{#1})}
\newcommand{\nc}{\newcommand}
\newcommand{\rnc}{\renewcommand}
\nc{\cal}{\mathcal}
\nc{\goth}{\mathfrak}
\rnc{\bold}{\mathbf}
\renewcommand{\frak}{\mathfrak}
\renewcommand{\Bbb}{\mathbb}

\nc{\Cal}{\cal}
\nc{\Xp}[1]{X^+(#1)}
\nc{\Xm}[1]{X^-(#1)}
\nc{\on}{\operatorname}
\nc{\ch}{\mbox{ch}}
\nc{\Z}{{\bold Z}}
\nc{\J}{{\cal J}}
\nc{\C}{{\bold C}}
\nc{\Q}{{\bold Q}}
\renewcommand{\P}{{\cal P}}
\nc{\N}{{\Bbb N}}
\nc\beq{\begin{equation}}
\nc\enq{\end{equation}}
\nc\lan{\langle}
\nc\ran{\rangle}
\nc\bsl{\backslash}
\nc\mto{\mapsto}
\nc\lra{\leftrightarrow}
\nc\hra{\hookrightarrow}
\nc\sm{\smallmatrix}
\nc\esm{\endsmallmatrix}
\nc\sub{\subset}
\nc\ti{\tilde}
\nc\nl{\newline}
\nc\fra{\frac}
\nc\und{\underline}
\nc\ov{\overline}
\nc\ot{\otimes}
\nc\bbq{\bar{\bq}_l}
\nc\bcc{\thickfracwithdelims[]\thickness0}
\nc\ad{\text{\rm ad}}
\nc\Ad{\text{\rm Ad}}
\nc\Hom{\text{\rm Hom}}
\nc\End{\text{\rm End}}
\nc\Ind{\text{\rm Ind}}
\nc\Res{\text{\rm Res}}
\nc\Ker{\text{\rm Ker}}
\rnc\Im{\text{Im}}
\nc\sgn{\text{\rm sgn}}
\nc\tr{\text{\rm tr}}
\nc\Tr{\text{\rm Tr}}
\nc\supp{\text{\rm supp}}
\nc\card{\text{\rm card}}
\nc\bst{{}^\bigstar\!}
\nc\he{\heartsuit}
\nc\clu{\clubsuit}
\nc\spa{\spadesuit}
\nc\di{\diamond}

\nc\al{\alpha}
\nc\bet{\beta}
\nc\ga{\gamma}
\nc\de{\delta}
\nc\ep{\epsilon}
\nc\io{\iota}
\nc\om{\omega}
\nc\si{\sigma}
\rnc\th{\theta}
\nc\ka{\kappa}
\nc\la{\lambda}
\nc\ze{\zeta}

\nc\vp{\varpi}
\nc\vt{\vartheta}
\nc\vr{\varrho}

\nc\Ga{\Gamma}
\nc\De{\Delta}
\nc\Om{\Omega}
\nc\Si{\Sigma}
\nc\Th{\Theta}
\nc\La{\Lambda}

\nc\boa{\bold a}
\nc\bob{\bold b}
\nc\boc{\bold c}
\nc\bod{\bold d}
\nc\boe{\bold e}
\nc\bof{\bold f}
\nc\bog{\bold g}
\nc\boh{\bold h}
\nc\boi{\bold i}
\nc\boj{\bold j}
\nc\bok{\bold k}
\nc\bol{\bold l}
\nc\bom{\bold m}
\nc\bon{\bold n}
\nc\boo{\bold o}
\nc\bop{\bold p}
\nc\boq{\bold q}
\nc\bor{\bold r}
\nc\bos{\bold s}
\nc\bou{\bold u}
\nc\bov{\bold v}
\nc\bow{\bold w}
\nc\boz{\bold z}

\nc\ba{\bold A}
\nc\bb{\bold B}
\nc\bc{\bold C}
\nc\bd{\bold D}
\nc\be{\bold E}
\nc\bg{\bold G}
\nc\bh{\bold H}
\nc\bi{\bold I}
\nc\bj{\bold J}
\nc\bk{\bold K}
\nc\bl{\bold L}
\nc\bm{\bold M}
\nc\bn{\bold N}
\nc\bo{\bold O}
\nc\bp{\bold P}
\nc\bq{\bold Q}
\nc\br{\bold R}
\nc\bs{\bold S}
\nc\bt{\bold T}
\nc\bu{\bold U}
\nc\bv{\bold V}
\nc\bw{\bold W}
\nc\bz{\bold Z}
\nc\bx{\bold X}

\nc\ca{\Cal A}
\nc\cb{\Cal B}
\nc\cc{\Cal C}
\nc\cd{\Cal D}
\nc\ce{\Cal E}
\nc\cf{\Cal F}
\nc\cg{\Cal G}
\rnc\ch{\Cal H}
\nc\ci{\Cal I}
\nc\cj{\Cal J}
\nc\ck{\Cal K}
\nc\cl{\Cal L}
\nc\cm{\Cal M}
\nc\cn{\Cal N}
\nc\co{\Cal O}
\nc\cp{\Cal P}
\nc\cq{\Cal Q}
\nc\car{\Cal R}
\nc\cs{\Cal S}
\nc\ct{\Cal T}
\nc\cu{\Cal U}
\nc\cv{\Cal V}
\nc\cz{\Cal Z}
\nc\cx{\Cal X}
\nc\cy{\Cal Y}

\nc\e[1]{E_{#1}}
\nc\ei[1]{E_{\delta - \alpha_{#1}}}
\nc\esi[1]{E_{s \delta - \alpha_{#1}}}
\nc\eri[1]{E_{r \delta - \alpha_{#1}}}
\nc\ed[2][]{E_{#1 \delta,#2}}
\nc\ekd[1]{E_{k \delta,#1}}
\nc\emd[1]{E_{m \delta,#1}}
\nc\erd[1]{E_{r \delta,#1}}

\nc\ef[1]{F_{#1}}
\nc\efi[1]{F_{\delta - \alpha_{#1}}}
\nc\efsi[1]{F_{s \delta - \alpha_{#1}}}
\nc\efri[1]{F_{r \delta - \alpha_{#1}}}
\nc\efd[2][]{F_{#1 \delta,#2}}
\nc\efkd[1]{F_{k \delta,#1}}
\nc\efmd[1]{F_{m \delta,#1}}
\nc\efrd[1]{F_{r \delta,#1}}

\nc\fa{\frak a}
\nc\fb{\frak b}
\nc\fc{\frak c}
\nc\fd{\frak d}
\nc\fe{\frak e}
\nc\ff{\frak f}
\nc\fg{\frak g}
\nc\fh{\frak h}
\nc\fj{\frak j}
\nc\fk{\frak k}
\nc\fl{\frak l}
\nc\fm{\frak m}
\nc\fn{\frak n}
\nc\fo{\frak o}
\nc\fp{\frak p}
\nc\fq{\frak q}
\nc\fr{\frak r}
\nc\fs{\frak s}
\nc\ft{\frak t}
\nc\fu{\frak u}
\nc\fv{\frak v}
\nc\fz{\frak z}
\nc\fx{\frak x}
\nc\fy{\frak y}

\nc\fA{\frak A}
\nc\fB{\frak B}
\nc\fC{\frak C}
\nc\fD{\frak D}
\nc\fE{\frak E}
\nc\fF{\frak F}
\nc\fG{\frak G}
\nc\fH{\frak H}
\nc\fJ{\frak J}
\nc\fK{\frak K}
\nc\fL{\frak L}
\nc\fM{\frak M}
\nc\fN{\frak N}
\nc\fO{\frak O}
\nc\fP{\frak P}
\nc\fQ{\frak Q}
\nc\fR{\frak R}
\nc\fS{\frak S}
\nc\fT{\frak T}
\nc\fU{\frak U}
\nc\fV{\frak V}
\nc\fZ{\frak Z}
\nc\fX{\frak X}
\nc\fY{\frak Y}
\nc\tfi{\ti{\Phi}}
\nc\bF{\bold F}
\rnc\bol{\bold 1}

\nc\ua{{_\ca\bold U}}

\nc\qinti[1]{[#1]_i}
\nc\q[1]{[#1]_q}
\nc\xpm[2]{E_{#2 \delta \pm \alpha_#1}}  
\nc\xmp[2]{E_{#2 \delta \mp \alpha_#1}}
\nc\xp[2]{E_{#2 \delta + \alpha_{#1}}}
\nc\xm[2]{E_{#2 \delta - \alpha_{#1}}}
\nc\hik{\ed{k}{i}}
\nc\hjl{\ed{l}{j}}
\nc\qcoeff[3]{\left[ \begin{smallmatrix} {#1}& \\ {#2}& \end{smallmatrix}
\negthickspace \right]_{#3}}
\nc\qi{q}
\nc\qj{q}

\nc\ufdm{{_\ca\bu}_{\rm fd}^{\le 0}}


\nc\isom{\cong} 

\nc{\pone}{{\Bbb C}{\Bbb P}^1}
\nc{\pa}{\partial}
\def\H{\cal H}
\def\L{\cal L}
\nc{\F}{{\cal F}}
\nc{\Sym}{{\goth S}}
\nc{\A}{{\cal A}}
\nc{\arr}{\rightarrow}
\nc{\larr}{\longrightarrow}

\nc{\ri}{\rangle}
\nc{\lef}{\langle}
\nc{\W}{{\cal W}}
\nc{\uqatwoatone}{{U_{q,1}}(\su)}
\nc{\uqtwo}{U_q(\goth{sl}_2)}
\nc{\dij}{\delta_{ij}}
\nc{\divei}{E_{\alpha_i}^{(n)}}
\nc{\divfi}{F_{\alpha_i}^{(n)}}
\nc{\Lzero}{\Lambda_0}
\nc{\Lone}{\Lambda_1}
\nc{\ve}{\varepsilon}
\nc{\phioneminusi}{\Phi^{(1-i,i)}}
\nc{\phioneminusistar}{\Phi^{* (1-i,i)}}
\nc{\phii}{\Phi^{(i,1-i)}}
\nc{\Li}{\Lambda_i}
\nc{\Loneminusi}{\Lambda_{1-i}}
\nc{\vtimesz}{v_\ve \otimes z^m}

\nc{\asltwo}{\widehat{\goth{sl}_2}}
\nc\ag{\widehat{\goth{g}}}  
\nc\teb{\tilde E_\boc}
\nc\tebp{\tilde E_{\boc'}}

\title[Monomial Bases of Quantized Enveloping Algebras]
{Monomial Bases of Quantized Enveloping Algebras}

\author{ Vyjayanthi Chari and Nanhua Xi }
\address{V.C. and N.X. Department of Mathematics, University of California,
Riverside, CA 92521, USA}
\address{N.X. Institute of Mathematics, Chinese Academy of Sciences,
Beijing 100080, China}
\email{chari@math.ucr.edu, nanhua@math.ucr.edu, nanhua@math08.math.ac.cn}





\dedicatory{}


\begin{abstract}
We construct a monomial basis of the positive part $\bu^+$ of the quantized
enveloping algebra associated to a finite--dimensional simple Lie algebra. As an 
application  we give a simple
proof of the existence and uniqueness of the canonical basis of $\bu^+$.
\end{abstract}

\maketitle

\section{\bf Introduction}

In [L1], Lusztig showed that the positive part $\bu^+$ of the quantized
enveloping algebra associated to a finite--dimensional simple Lie algebra 
$\frak{g}$, had a remarkable basis called the canonical basis. The main idea in 
proving its existence and uniqueness  was the following.
Corresponding to every reduced expression $\bold{i}$ of the longest element 
$w_0$ of the
Weyl group of $\frak{g}$ one constructs a Poincare--Birkhoff--Witt basis 
$B_{\bold i}$ of $\bu^+$.  Lusztig proved in [L1] that the 
$\bz[q^{-1}]$--lattice $\cl$ spanned by $B_i$ is independent of the  choice of 
$\bold{i}$  and that the image of $B_{\bold i}$ in the $\bz$--module 
$\cal{L}/q^{-1}L$ is a basis $B$  independent of $\bold{i}$. Let 
$\overline{\cal{L}}$ be the image of $\cal{L}$ under the bar map  (a certain 
$\bq$--algebra involution) of $\bu^+$. The canonical basis $\bold{B}$ is the 
preimage of $B$ in $\cal{L}\cap\overline{\cal{L}}$.

In [K1], Kashiwara introduced the notion of crystal bases for the quantized 
algebras of classical type and later generalized it to quantized algebras 
associated to an arbitrary symmetrizable Kac--Moody algebra. The crystal basis 
is a base {\hbox{\lq\ at $q=0$'} of $\bu^+$ with certain properties. Later, he  
proved that the crystal bases could be  \lq melted\rq{} to give a basis of 
$\bu^+$ itself, which is called the global crystal base. The main tool used here 
was a certain bilinear form on the algebra, and the global crystal basis can be 
characterized as   a bar--invariant quasi--orthonormal basis with respect to 
this form. In [X2], Xi proved that the bases $B_{\bold i}$
are quasi--orthonormal with respect to this form.

The quantized enveloping algebra also admits another  symmetric bilinear form 
introduced by Drinfeld, and  
Lusztig proved in [L2] that the bases $B_{\bold i}$ are quasi--orthonormal and 
that   the canonical basis can be characterized as the  bar--invariant, 
quasi--orthonormal basis of $\bu^+$ with respect to the Drinfeld form.  It is 
now proved [GL] that in fact the global crystal basis and the canonical basis 
are the same.

In this paper we construct  a basis of $\bu^+$ whose terms are monomials in the 
Chevalley generators and hence  is bar--invariant but not quasi--orthonormal.  
We are then able to give a very simple proof of the existence and uniqueness of 
the canonical basis.
Our construction of the monomial basis depends on picking a specific reduced 
expression for $w_0$. We conjecture that in fact there exists a monomial basis
corresponding to  every reduced expression. In view of [BCP] we expect also
that similar results should be true for the quantized affine algebras.

\medskip
\def\bqv{\bold Q(v)}
\def\a{\alpha}

\section{Preliminaries} Throughout this paper $\frak{g}$ will denote a 
finite-dimensional complex simple Lie algebra, $(a_{ij})_{i,j\in I},
\ I=\{1,\dots, n\},$ will denote its Cartan matrix and $\bn[I]$ will denote the 
set of linear combinations $\nu=\sum_i\nu_ii$, $\nu_i\in\bz, \nu_i\ge 0$. Let 
$W$ be the Weyl group of $\frak{g}$. It is well--known that $W$ is a Coxeter 
group generated by simple reflections $s_i$ for $i\in I$. Let $l(w)$ denote the 
length of a reduced expression of $w$ and let $w_0$ be the unique element of 
maximal length in $W$.
Let $R^+$ denote a  set of positive roots of $\frak{g}$ and let $\alpha_1,\dots, 
\alpha_n$ be the  set of simple roots.
Corresponding to any reduced expression of $w_0$, say $w_0=s_{k_1}s_{k_2}\dots 
s_{k_r}$, $r=|R^+|$,  we have a total order on $R^+$,
\begin{equation*} \beta_1<\beta_2\cdots <\beta_r\end{equation*}
where $\beta_i=s_{k_1}s_{k_2}\dots s_{k_{i-1}}\alpha_{k_i}$.

For the rest of the paper we shall be  working with 
two specific reduced expressions of $w_0$.  Since the numbering of the nodes is 
important for our purposes, we include below,  for the readers convenience, the 
Dynkin diagrams of the  various finite--dimensional complex simple Lie algebras. 
We assume that the node $\alpha_1$ is short if $\frak{g}$ is of type
 $F_4$ or $G_2$.

\def\no{\noindent}

 Type $A_n\ (n\geq 1)$. 

\

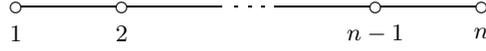
\begin{picture}(200,30)(-60,0)
\multiput(10,20)(40,0){2}{\circle{4}}
\multiput(12,20)(40,0){2}{\line(1,0){36}}
\multiput(93,20)(5,0){3}{\line(1,0){1}}
\multiput(108,20)(40,0){2}{\line(1,0){36}}
\multiput(146,20)(40,0){2}{\circle{4}}
\put(10,10){\makebox(0,0)[c]{\small 1}}
\put(50,10){\makebox(0,0)[c]{\small 2}}
\put(146,10){\makebox(0,0)[c]{{\small $n-1$}}}
\put(186,10){\makebox(0,0)[c]{{\small $n$}}}
\end{picture}

 Type $B_n\ (n \geq 2).$ 

\

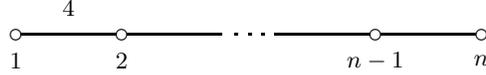
\begin{picture}(180,30)(-60,0)
\multiput(10,20)(40,0){2}{\circle{4}}
\multiput(12,20)(40,0){2}{\line(1,0){36}}
\multiput(93,20)(5,0){3}{\line(1,0){1}}
\multiput(108,20)(40,0){2}{\line(1,0){36}}
\multiput(146,20)(40,0){2}{\circle{4}}
\put(10,10){\makebox(0,0)[c]{\small 1}}
\put(50,10){\makebox(0,0)[c]{\small 2}}
\put(146,10){\makebox(0,0)[c]{{\small $n-1$}}}
\put(186,10){\makebox(0,0)[c]{{\small $n$}}}
\put(30,30){\makebox(0,0)[c]{\small 4}}
\end{picture}

 Type $C_n\ (n \geq 3).$ 

\

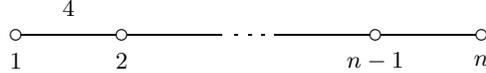
\begin{picture}(180,30)(-60,0)
\multiput(10,20)(40,0){2}{\circle{4}}
\multiput(12,20)(40,0){2}{\line(1,0){36}}
\multiput(93,20)(5,0){3}{\line(1,0){1}}
\multiput(108,20)(40,0){2}{\line(1,0){36}}
\multiput(146,20)(40,0){2}{\circle{4}}
\put(10,10){\makebox(0,0)[c]{\small 1}}
\put(50,10){\makebox(0,0)[c]{\small 2}}
\put(146,10){\makebox(0,0)[c]{{\small $n-1$}}}
\put(186,10){\makebox(0,0)[c]{{\small $n$}}}
\put(30,30){\makebox(0,0)[c]{\small 4}}
\end{picture}\\

 Type $D_n\ (n \geq 4).$ 

\

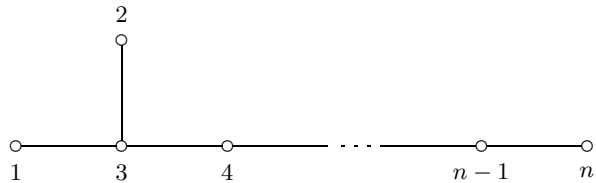
\begin{picture}(240,90)(-60,0)
\multiput(10,20)(40,0){3}{\circle{4}}
\multiput(12,20)(40,0){3}{\line(1,0){36}}
\multiput(133,20)(5,0){3}{\line(1,0){1}}
\put(50,60){\circle{4}}
\put(50,22){\line(0,1){36}}
\put(50,70){\makebox(0,0)[c]{\small 2}}
\multiput(148,20)(40,0){2}{\line(1,0){36}}
\multiput(186,20)(40,0){2}{\circle{4}}
\put(10,10){\makebox(0,0)[c]{\small 1}}
\put(50,10){\makebox(0,0)[c]{\small 3}}
\put(90,10){\makebox(0,0)[c]{\small 4}}
\put(186,10){\makebox(0,0)[c]{{\small $n-1$}}}
\put(226,10){\makebox(0,0)[c]{{\small $n$}}}
\end{picture}\\

 Type $F_4$

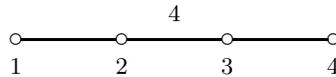
\begin{picture}(180,30)(-60,0)
\multiput(10,20)(40,0){4}{\circle{4}}
\multiput(12,20)(40,0){3}{\line(1,0){36}}
\put(70,30){\makebox(0,0)[c]{\small 4}}
\put(10,10){\makebox(0,0)[c]{\small 1}}
\put(50,10){\makebox(0,0)[c]{\small 2}}
\put(90,10){\makebox(0,0)[c]{\small 3}}
\put(130,10){\makebox(0,0)[c]{\small 4}}
\end{picture}\\

\vfill\eject

 Type $G_2$

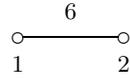
\begin{picture}(180,30)(-60,0)
\multiput(10,20)(40,0){2}{\circle{4}}
\put(12,20){\line(1,0){36}}
\put(30,30){\makebox(0,0)[c]{\small 6}}
\put(10,10){\makebox(0,0)[c]{\small 1}}
\put(50,10){\makebox(0,0)[c]{\small 2}}
\end{picture}

\

 Type $E_6 $

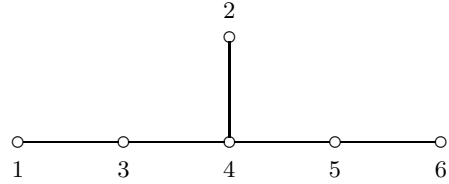
\begin{picture}(180,90)(-60,0)
\multiput(10,20)(40,0){5}{\circle{4}}
\multiput(12,20)(40,0){4}{\line(1,0){36}}
\put(90,60){\circle{4}}
\put(90,22){\line(0,1){36}}
\put(10,10){\makebox(0,0)[c]{\small 1}}
\put(50,10){\makebox(0,0)[c]{\small 3}}
\put(90,10){\makebox(0,0)[c]{\small 4}}
\put(130,10){\makebox(0,0)[c]{\small 5}}
\put(170,10){\makebox(0,0)[c]{\small 6}}
\put(90,70){\makebox(0,0)[c]{\small 2}}
\end{picture}

\

 Type $E_7$ 

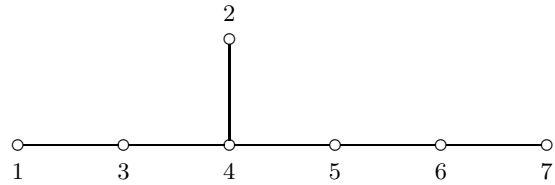
\begin{picture}(220,90)(-60,0)
\multiput(10,20)(40,0){6}{\circle{4}}
\multiput(12,20)(40,0){5}{\line(1,0){36}}
\put(90,60){\circle{4}}
\put(90,22){\line(0,1){36}}
\put(10,10){\makebox(0,0)[c]{\small 1}}
\put(50,10){\makebox(0,0)[c]{\small 3}}
\put(90,10){\makebox(0,0)[c]{\small 4}}
\put(130,10){\makebox(0,0)[c]{\small 5}}
\put(170,10){\makebox(0,0)[c]{\small 6}}
\put(210,10){\makebox(0,0)[c]{\small 7}}
\put(90,70){\makebox(0,0)[c]{\small 2}}
\end{picture}

\

 Type $E_8$ 

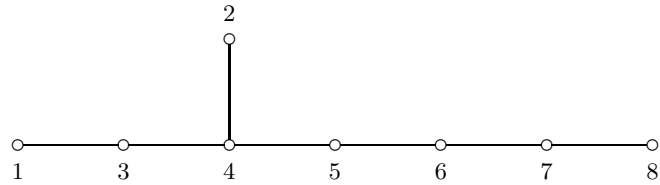
\begin{picture}(260,90)(-60,0)
\multiput(10,20)(40,0){7}{\circle{4}}
\multiput(12,20)(40,0){6}{\line(1,0){36}}
\put(90,60){\circle{4}}
\put(90,22){\line(0,1){36}}
\put(10,10){\makebox(0,0)[c]{\small 1}}
\put(50,10){\makebox(0,0)[c]{\small 3}}
\put(90,10){\makebox(0,0)[c]{\small 4}}
\put(130,10){\makebox(0,0)[c]{\small 5}}
\put(170,10){\makebox(0,0)[c]{\small 6}}
\put(210,10){\makebox(0,0)[c]{\small 7}}
\put(250,10){\makebox(0,0)[c]{\small 8}}
\put(90,70){\makebox(0,0)[c]{\small 2}}
\end{picture}
\medskip

Let $\gamma_n\in W$ be defined as follows:
\begin{equation*} \gamma_n= \begin{cases}s_1s_{2}\cdots s_n & \text{if 
$\frak{g}$ is of type $A_n$,}\\
s_ns_{n-1}\cdots s_2s_1s_2s_3\cdots s_{n-1}s_{n} &\text{if $\frak{g}$ is of type 
$B_n$ or $C_n$,}\\
s_ns_{n-1}\cdots s_3s_1s_2s_3\cdots s_{n-1}s_n &\text{if $\frak{g}$ is of type 
$D_n$.}\end{cases}
\end{equation*}
It is easy to see that  \begin{equation*}\bold{j}=\gamma_n\gamma_{n-1}\cdots 
\gamma_1\end{equation*} is a reduced expression of $w_0$.
 
For the exceptional algebras, we take $\bold{j}$
 as below:

\begin{eqnarray*}
G_2& s_2s_1s_2s_1s_2s_1,\\
F_4& s_4s_3s_2s_3s_4s_1s_2s_3s_2s_1s_4s_3s_2s_3s_4s_1s_2s_3s_2s_1s_2s_3s_2s_3,
\end{eqnarray*}
\begin{eqnarray*}
E_6& s_1s_3s_4s_2s_5s_4s_3s_1s_6s_5s_4s_2s_3s_4s_5s_6u_1,\\
E_7&s_7s_6s_5s_4s_2s_3s_4s_5s_6s_7s_1s_3s_4s_5\\&\times 
s_2s_4s_3s_1s_6s_5s_4s_2s_3s_4s_5s_6s_7u_2,\\
E_8&
s_8s_7s_6s_5s_4s_2s_3s_4s_5s_6s_7s_8s_1s_3s_4s_2s_5s_4s_3s_1s_6s_5s_7s_6s_4s_3s_
2s_5s_4s_5s_3s_2s_4\\&\times 
s_6s_5s_7s_6s_1s_3s_4s_2s_5s_4s_3s_1s_8s_7s_6s_5s_4s_2s_3s_4s_5s_6s_7s_8u_3,
\end{eqnarray*}
where $u_1$ (resp. $u_2,u_3$)  is the reduced expression of the longest element 
of $D_5$ (resp. $E_6, E_7$ ) obtained by dropping the node 6 (resp. 7,8) which 
has been chosen previously.

If $\frak{g}$
 is of type $A_n$, $D_n$ or $E_6$, let $\tau$ be the non--trivial diagram 
automorphism of order 2 and let $\bold{i}$ be the reduced expression for $w_0$ 
obtained by applying $\tau$ to $\bold{j}$; for $\frak{g}$ of other types, we 
take $\bold{j}=\bold{i}$.
In what follows we assume that the roots $\beta_i$, $i=1,\dots , r$ are defined
with respect to the reduced expression $\bold{i}$.

Let $q$ be an indeterminate, let $\bq(q)$ be the field of rational
functions in $q$ with rational coefficients, and let $\bz[q,q^{-1}]$ be
the ring of Laurent polynomials with integer coefficients. For
$r,m,d\in\bn$, $m\ge r$, define
\begin{equation*}
 [m]_d=\frac{q^{dm} -q^{-dm}}{q^d -q^{-d}},\ \ \ \ 
  [m]_d! =[m]_d[m-1]_d\ldots [2]_d[1]_d,\ \ \ \ 
\left[\begin{matrix} m\\ 
  r\end{matrix}\right]_d = \frac{[m]_d!}{[r]_d![m-r]_d!}.
\end{equation*}
  Then $\left[\begin{matrix} m\\r\end{matrix}\right]_d\in\bz[q,q^{-1}]$
  for all $m\ge r\ge 0$. Choose $d_i\in\{1,2,3\}$ such that $(d_ia_{ij})$ is 
symmetric and such that $\sum_i d_i$ is minimal.

\begin{prop}{\label{defnbu}} There is a Hopf algebra $\bu$ over $\bq(q)$ which 
is generated as an algebra by elements $E_{\alpha_i}$, $F_{\alpha_i}$, 
$K_i^{{}\pm 1}$ ($i\in I$), with the following defining relations:
\begin{align*} 
  K_iK_i^{-1}=K_i^{-1}K_i&=1,\ \ \ \ K_iK_j=K_jK_i,\\ 
  K_iE_{\alpha_j} K_i^{-1}&=q^{ d_ia_{ij}}E_{\alpha_j},\\ 
K_iF_{\alpha_j} K_i^{-1}&=q^{-d_ia_{ij}}F_{\alpha_j},\\
  [E_{\alpha_i}, F_{\alpha_j}
]&=\delta_{ij}\frac{K_i-K_i^{-1}}{q^{d_i}-q^{-d_i}},\\ 
  \sum_{r=0}^{1-a_{ij}}(-1)^r\left[\begin{matrix} 1-a_{ij}\\ 
  r\end{matrix}\right]_{d_i}
&(E_{\alpha_i})^rE_{\alpha_j}(E_{\alpha_i})^{1-a_{ij}-r}=0\ 
  \ \ \ \ \text{if $i\ne j$},\\
\sum_{r=0}^{1-a_{ij}}(-1)^r\left[\begin{matrix} 1-a_{ij}\\ 
  r\end{matrix}\right]_{d_i}
&(F_{\alpha_i})^rF_{\alpha_j}(F_{\alpha_i})^{1-a_{ij}-r}=0\ 
  \ \ \ \ \text{if $i\ne j$}.
\end{align*}
The comultiplication of $\bu$ is given on generators by
$$\Delta(E_{\alpha_i})=E_{\alpha_i}\ot 1+K_i\ot E_{\alpha_i},\ \ 
\Delta(F_{\alpha_i})=F_{\alpha_i}\ot K_i^{-1} + 1\ot F_{\alpha_i},\ \ 
\Delta(K_i)=K_i\ot K_i,$$
for $i\in I$.\hfill\qedsymbol
\end{prop}

Let $\bu^+$   be the $\bq(q)$-subalgebra of $\bu$ generated by the 
$E_{\alpha_i}$  for $i\in {I}$.

\begin{defn} An element $x\in\bu^+$ is called bar--invariant if   it is fixed by 
the  $\bq$--algebra homomorphism $^-:\bu^+\to\bu^+$ defined by extending:
\begin{equation*} \overline{E_{\alpha_i}} = E_{\alpha_i},\ \ \  \overline{q} 
=q^{-1}.\end{equation*}
\end{defn}

For $\nu\in\bn[I]$, let $\bu^+_\nu$ be the subspace of $\bu^+$ spanned by the 
monomials \break $E_{\alpha_{s_1}}E_{\alpha_{s_2}}\cdots E_{\alpha_{s_t}}$ such that 
for any $i\in I$, the number of occurrences of $i$ in the sequence $s_1,\dots, 
s_t$ is equal to $\nu_i$. An element $x\in\bu^+$ is said to have homogeneity 
$\nu$ if $x\in\bu^+_\nu$ and we denote its homogeinty by $|x|$.

It is convenient to use the following notation:
\begin{equation*}
E_{\alpha_i}^{(r)}=\frac{E_{\alpha_i}^r}{[r]_{d_i}!}.
\end{equation*} 
The elements $F_{\alpha_i}^{(r)}$ are defined similarly. 

Set $\ca= Z[q,q^{-1}]$ and let $\ua^+$ be the $\ca$--subalgebra of $\bu^+$ 
generated by $E_{\alpha_i}^{(r)}$, $i\in I$, $r\ge 0$.  

For $i\in I$, let $T_i$ ($i\in I$) be
  the $\bq(q)$-algebra automorphisms of $\bu$ defined as follows (see [L2]):
  \begin{align*}\label{braid}
    T_i(E_{\alpha_i}^{(m)})&=(-1)^mq^{-m(m-1)}K_i^{-m}F_{\alpha_i}^{(m)},\\ 
    T_i(F_{\alpha_i}^{(m)})&=(-1)^mq^{m(m-1)}E_{\alpha_i}^{(m)}K_i^{m},\\ 
    T_i(E_{\alpha_j}^{(m)})&=\sum_{r=0}^{-ma_{ij}}(-1)^{r}q^{-r}
    E_{\alpha_i}^{(r)}E_{\alpha_j}^{(m)}E_{\alpha_i}^{(-ma_{ij}-r)}\ 
    \ \text{if $i\ne j$},\\ 
    T_i(F_{\alpha_j}^{(m)})&=\sum_{r=0}^{-ma_{ij}}(-1)^{r}q^{r}
    F_{\alpha_i}^{(-ma_{ij}-r)}
    F_{\alpha_j}^{(m)}F_{\alpha_i}^{(r)}\ \ {\text{if
        $i\ne j$}}.
\end{align*}
Recall that $\bold i=s_{i_1}s_{i_2}\cdots s_{i_r}$.
For $j=1,2,\dots, r$, $m\in\bz$, $m\ge 0$, define a set of root vectors by 
\begin{equation*}
E_{\beta_j}^{(m)} =T_{i_1}T_{i_2}\dots 
T_{i_{j-1}}(E_{\alpha_{i_j}}^{(m)}).\end{equation*}

For $\boc=(c_1,...,c_r)\in \bn^r$, we set $
E_\boc=
E_{\beta_1}^{(c_1)}E_{\beta_2}^{(c_2)}\cdots E_{\beta_r}^{(c_r)}$. 
The elements $E_\boc$ are clearly homogenous and $|E_\boc| 
=\sum_{j=1}^rc_j\beta_j$. 
The following result is proved in [L3, Corollary 40.2.2]. 
\begin{prop}{\label{PBW}} The set $\{E_{\boc}:\boc\in \bn^r\}$ is  a 
$\bq(q)$--basis of $\bu^+$ and an $\ca$--basis of $\ua^+$.
\qedsymbol\end{prop}
\vskip 12pt

\section{ Monomial basis}
\def\ds{\displaystyle\sum}

In this section we construct a monomial basis of $\bu^+$, i.e, a basis 
consisting
of products of $E_{\alpha_i}^{(s)}$, $i\in I$, $s\ge 0$. We  also give
a simple proof for the existence and uniqueness of the canonical basis
of $\bu^+$. Recall that we have fixed a reduced expression $\bold{j}= 
s_{j_1}s_{j_2}\cdots s_{j_r}$ of $w_0$. For $\boc=(c_1,...,c_r)\in
\bn^r$ we set \begin{equation*} M_\boc=E_{j_1}^{(c_1)}E_{j_2}^{(c_2)}
\cdots E_{j_r}^{(c_r)}.\end{equation*} Let $>$ be the  lexicographic ordering on 
$\bn^r$ such that 
\begin{equation*}(1,0,\cdots ,0)> (0,1,0,\cdots ,0)>\cdots >(0,0,\cdots 
,1).\end{equation*}
\begin{thm}{\label{mbases}}
\begin{enumerate}

\item[(i)] 
For any $\boc\in\bn^r$ there exists
$f(\boc)\in\bn^r$ such that \begin{equation*} M_{f(\boc)}=E_\boc+\ds_{
\stackrel{\scriptstyle \bod\in\bn^r} {\bod>\boc}}  \xi_\bod E_\bod 
,\end{equation*} where $\xi_\bod\in\A$. 
\item[(ii)] The set $\{M_{f(\boc)}: \boc\in\bn^r\}$ is an $\A$--basis of 
$\ua^+$.
\end{enumerate}

\end{thm}

We postpone the proof of the theorem and deduce first the result on canonical 
bases.
Let $\cal L$ be the $\bold Z[q^{-1}]$--lattice of $\bold U^+$ spanned by the set 
$\{E_{\boc}\ :\ \boc\in\bn^r\}$.
\begin{thm}[L3,K2] For each $\boc\in\bn^r$  there exist a unique $b_\boc$ in the 
lattice $ \cl$  such that
$b_\boc$ is bar--invariant and 
\begin{equation*}b_\boc=E_\boc+\ds_{\stackrel{\scriptstyle \bod\in\bn^r
}{\bod>\boc}}\zeta_\bod E_\bod ,\end{equation*}
where $\zeta_\bod\in q^{-1}\bold Z[q^{-1}]$. The set
\begin{equation*}\bold{B}=\{b_\boc:\boc\in\bn^r\}\end{equation*} is an 
$\ca$--basis of $\ua^+$ and  is called the canonical basis or the global crystal 
basis of $\bu^+$. 
\end{thm}
\begin{pf} For each $\boc\in\bn^r$, observe that the set
\begin{equation*} \bold{S_\boc} =\{\bod\in \bn^r:\bod\ge \boc, \ \ |E_\boc| 
=|E_\bod|\}\end{equation*}
is finite and not empty. Let  $\boc =\boc_0<\boc_1<\cdots<\boc_m$ be the 
elements of $S_\boc$. If $\boc_0=\boc_m$ 
then by Theorem 1 
\begin{equation*} 
E_\boc =M_{f(\boc)}\end{equation*}
and there is nothing to prove. Otherwise, by Theorem 1
we can write
\begin{equation*} M_{f(\boc)}=E_\boc+\sum_{1\le k\le m} 
\xi_kE_{\boc_k},\end{equation*}
for some $\xi_k\in Z[q,q^{-1}]$. Let $\xi_1'$ be the unique bar--invariant 
element of  $\ca$ such that
\begin{equation*} \xi_1-\xi_1'=\eta_1\in q^{-1}Z[q^{-1}].\end{equation*}
Applying Theorem 1, we get
\begin{equation*} M_{f(\boc)}-\xi_1'M_{f(\boc_1)} 
=E_{\boc}+\eta_1E_{\boc_1}+\sum_{2\le k\le m}\xi_{k,2}E_{\boc_k},\end{equation*}
where $\xi_{k,2}\in \bz[q,q^{-1}]$ for all $2\le k\le m$. Next, let  $\xi_2'$ be 
the unique bar--invariant element of $\bz[q,q^{-1}]$  such that 
\begin{equation*}\xi_{2,2}-\xi_2'=\eta_2\in q^{-1}\bz[q^{-1}].\end{equation*} As 
before, we can find elements $\xi_{k,3}\in\bz[q,q^{-1}]$ for $3\le k\le m$ such 
that,
\begin{equation*}M_{f(\boc)}-\xi_1'M_{f(\boc_1)} -\xi_2'M_{f(\boc_2)} = 
E_{\boc}+\eta_1E_{\boc_1}+\eta_2E_{\boc_2}+\sum_{3\le k\le 
m}\xi_{k,3}E_{\boc_k}. \end{equation*} 
Repeating this process we find finally that
\begin{equation*}  M_{f(\boc)} -\sum_{1\le k\le m} \xi_k'M_{f(\boc_k)} = 
E_\boc+\sum_{1\le k\le m} \eta_kE_{\boc_k},\end{equation*} with $\eta_k\in 
q^{-1}\bz[q^{-1}]$.  Since the left-hand side in the previous equation is 
obviously bar--invariant,
 the result follows by taking 
\begin{equation*} b_\boc= M_{f(\boc)} -\sum_{1\le k\le m} 
\xi_k'M_{f(\boc_k)}.\end{equation*}\end{pf}

We now turn to the proof of Theorem 1.
We assume (i) and prove (ii). Suppose that
\begin{equation*}
\sum_\boc \xi_{\boc} M_{f(\boc)} =0,
\end{equation*} 
for some 
$\xi_\boc\in\ca$.  Choose, if possible, $\boc_0\in\bn^r$ minimal such that 
$\xi_{\boc_0}\ne 0$. Using Theorem 1(i) we get that
\begin{equation*}
\xi_{\boc_0}E_{\boc_0}+\sum_{\boc>\boc_0}\xi_{\boc}E_{\boc}=0.
\end{equation*}
But this contradicts Proposition 1.2 and hence $\xi_{\boc_0} =0$  proving that 
the elements $M_{f(\boc)}\in\ua^+$ are linearly independent.

For each $\eta\in\bold N[I]$, set \begin{equation*}\ua^+_\eta 
=\ua^+\cap\bu^+_\eta. \end{equation*}
Let
\begin{equation*} S_\eta=\{\boc\in\bn^r: E_\boc\in\ua^+_\eta\},\end{equation*}
and let $\boc_0<\boc_1<\cdots<\boc_m$ be the elements of $S_\eta$.   From 
Theorem 1(i) it is clear that
\begin{equation*} M_{f(\boc_m)} =E_{\boc_m}.\end{equation*}
An obvious downward induction on $|S_\eta|$ proves that $E_{\boc_k}$ is in the 
span  of $M_{f(\bod)}$ for $\bod\in\bn^r$ and the result follows by Proposition 
1.2. 

\vskip 12pt

The rest of the paper is devoted to  proving Theorem 1(i). We need the following 
result proved in [L1, Section 6].  Let $\alpha\in R^+$ be such  that
\begin{equation*}\alpha-\alpha_n=\sum r_i\alpha_i,\ \ r_i\in \bz, \ \  r_i\geq 
0.\end{equation*}  Fix $k<n$ and let 
\begin{equation*} R^+_{\alpha, k} = \{\beta\in R^+: \beta =r\alpha+s\alpha_k, 
r,s\in\bz\}.\end{equation*}

\def\ds{\displaystyle\sum}

\medskip
\begin{lem}{\label{commute}} Let $a,b\in\bold N$.
\begin{enumerate}
\item[(i)] If $R^+_{\alpha,k}=\{\a,\a_k\}$, then 
\begin{equation*}
E_{\a_k}^{(a)}E_\a^{(b)}=E_\a^{(b)}E_{\a_k}^{(a)}.
\end{equation*}

\item[(ii)] If $R^+_{\alpha,k}=\{\a,\a+\a_k,\a_k\}$, then $d_\a=d_{\a_k}=d$ and 
 \begin{eqnarray*} E_{\a+\a_k}^{(a)} E^{(b)}_\a &= q^{dab} E^{(b)}_{\a} 
E^{(a)}_{\a+\a_k}, \\
 E^{(a)}_{\a_k} E^{(b)}_{\a+\a_k} &= q^{dab} E^{(b)}_{\a+\a_k} E^{(a)}_{\a_k},
 \\
 E^{(a)}_{\a_k} E^{(b)}_{\a} &= \ds_{r\in\bn}  q^{-d (a-r)(b-r)} E^{(b-r)}_{\a} 
E^{(r)}_{\a+\a_k} E^{(a-r)}_{\a_k}. 
 \end{eqnarray*}

\def\ds{\displaystyle\sum}
\item[(iii)] If $R^+_{\alpha,k}=\{\a,\a+\a_k,\a+2\a_k,\a_k\}$, then $d_\a=2$, 
$d_{\a_k}=1$ and 
\begin{eqnarray*}E^{(a)}_{\a+\a_k} E^{(b)}_{\a} &= q^{2ab} E^{(b)}_{\a} 
E^{(a)}_{\a+\a_k},\\
E^{(a)}_{\a+2\a_k} E^{(b)}_{\a+\a_k} &= q^{2ab} E^{(b)}_{\a+\a_k} 
E^{(a)}_{\a+2\a_k},\\
E^{(a)}_{\a_k} E^{(b)}_{\a+2\a_k} &= q^{2ab} E^{(b)}_{\a+2\a_k} E^{(a)}_{\a_k}, 
\end{eqnarray*} 
\begin{eqnarray*} E^{(a)}_{\a+2\a_k} E^{(b)}_{\a}& = \ds _{ r\in \bold N } q ^{2 
r(b-r)+ 2 r(a-r)}\displaystyle\prod ^ {r}_{h=1} 
(q^{4h -2} -1) E^{(b-r)}_{\a} E^{(2r)} _{\a+\a_k} E^{(a-r)}_{\a+2\a_k}, \\
 E^{(a)}_{\a_k} E^{(b)}_{\a+\a_k}& = \ds _{r\in \bold N } q^{r(b-r)+ r(a-r)- 
r}\displaystyle\prod ^ {r}_{h=1} 
(q^{2h}+1) E^{(b-r)}_{\a+\a_k} E^{(r)}_{\a+2\a_k} E^{(a-r)}_{\a_k}, \\
 E^{(a)}_{\a_k}E^{(b)}_\a &= \ds _{r,t\in \bn}
  q^{-2(b-r-t)(a-r-t)-(a-r-2t)r}E^{(b-r-t)}_{\alpha}
E^{(r)}_{\a+\a_k} E^{(t)}_{\a+2\a_k} E^{(a-r-2t)}_{\a_k}.\end{eqnarray*}
\end{enumerate}
 \hfill\qedsymbol\end{lem}

The proof of Theorem 1(i) proceeds by induction on the rank of $\frak{g}$. For 
rank one there is nothing to prove, and the rank two case is contained in the 
next Lemma. 
\begin{lem} 
\begin{enumerate}
\item[(i)] If $\frak{g}= A_2$, and $\boc=(c_1,c_2,c_3)$, then \begin{equation*} 
M_{f(\boc)} = E_1^{(c_2)}E_2^{(c_1+c_2)}E_1^{(c_3)}.\end{equation*}
\item[(ii)] If $\frak{g}=B_2$, and $\boc=(c_1,c_2,c_3,c_4)$, then,
\begin{equation*} M_{f(\boc)}= 
E_2^{(c_1)}E_1^{(c_2+2c_3)}E_2^{(c_2+c_3)}E_1^{(c_4)}.\end{equation*}
\item[(iii)] If $\frak{g}$ is of type $G_2$, and  $\boc=(c_1,c_2,\dots ,c_6)$,
then
\begin{equation*}M_{f(\boc)} 
=E_2^{(c_1)}E_1^{(c_2+3c_3)}E_2^{(c_2+2c_3)}E_1^{(2c_4+3c_5)}E_2^{(c_4+c_5)}E_1^
{(c_6)}.\end{equation*}
\end{enumerate}
\end{lem}

\begin{pf} If $\frak{g}$ is of type $A_2$ or $B_2$, the result follows by using 
the previous lemma. For $G_2$, see [X1].   \end{pf}

Turning to the general case, let $l=l_n$ be the number of positive roots $\beta$ 
such that $\beta-\alpha_n =\sum_{i=1}^n n_i\alpha_i$, where $n_i\ge 0$ for all 
$i=1,\ldots ,n$.  For any $\boc\in\bn^r$, we write 
\begin{equation*}\boc=\boc'+\boc'',\end{equation*}
where $c'_k=0$ if $k>l$ and $c''_k=0$ if $k\le l$.

\begin{lem}{\label{special}}
 Let $\boc\in\bn^r$ be such that $\boc''=0$.  Then there exists $f(\boc)$ such 
that \begin{equation*} M_{f(\boc)} = E_\boc+\sum_{\stackrel{\scriptstyle 
\bod\in\bn^r} {\bod'>\boc'}}\xi_{\bod} E_\bod\qquad 
\xi_{\bod}\in\bz[q,q^{-1}].\end{equation*}\end{lem}

Assuming the lemma we complete the proof as follows. 
 
 The element
$E_{\boc''}$ can be regarded as an element of $\bold{U}_q(\frak{g}_{n-1})$
where $\frak{g}_{n-1}$ is the simple Lie algebra associated to the 
$(n-1)\times (n-1)$ Cartan matrix obtained by dropping the $n^{th}$ row and the 
$n^{th}$ column. Hence by induction, we have a monomial $M_{f(\boc'')}$ such 
that
\begin{equation*} M_{f(\boc'')} = E_{\boc''}+\ds_{\bod>\boc''}\eta_\bod 
E_{\bod}, \end{equation*}
where $\eta_\bod\in\ca$ and $d_k=0$ if $k\le l$.
Let $M_{f(\boc')}$ be the monomial defined in the preceding lemma and set 
$$M_{f(\boc)} =M_{f(\boc')}M_{f(\boc'')}.$$

We get
\begin{eqnarray*} &M_{f(\boc)} =M_{f(\boc')}M_{f(\boc'')}\\
=&E_{\boc'}E_{\boc''}+E_{\boc'}\ds_{\bod>\boc''}\eta_\bod E_\bod
+\ds_{\stackrel{\scriptstyle \bod\in\bn^r} {\bod'>\boc'}}\xi_\bod E_\bod 
E_{\boc''}+\ds_{\stackrel{\scriptstyle \bod\in\bn^r} {\bod'>\boc'}}\xi_\bod 
E_\bod\ds_{\bod_1>\boc''}\eta_{\bod_1} E_{\bod_1}\\
=& E_{\boc}+\ds_{\stackrel{\scriptstyle \bod\in\bn^r}{\bod>\boc}}\zeta_\bod 
E_\bod,\qquad\ \ \ \zeta_\bod\in 
q^{-1}\bz[q^{-1}].\qquad\qedsymbol\end{eqnarray*}.

It remains to prove Lemma 2.3. Note that if $\boc''=0$, it is enough to prove 
that $\bod>\boc$ since this together with the fact that $|E_\bod|=|E_\boc|$ 
implies that $\bod'>\boc'$.
Let $j_1,j_2,\ldots, j_l$ be the first $l$--indices of the reduced expression 
$\bold{j}$ of $w_0$. Then using Lemma 2.1 repeatedly, it is not hard to see 
that, 
\begin{equation*} M_{f(\boc)} = E_{j_1}^{(k_1)} E_{j_2}^{(k_2)}\cdots 
E_{j_l}^{(k_l)},\end{equation*}
where $\bok=(k_1,k_2,\ldots, k_l)$ is related to $\boc$ as follows:

\begin{equation*}
k_i =c_n+c_{n-1}+\cdots+ c_{n-i+1}\end{equation*}
if $\frak{g}$ is of type $A_n$;
\begin{equation*} 
k_i = \begin{cases} c_i+c_{2n-i+1}+c_{2n-i+2}+\cdots+c_{2n-1}& {\text{if $2\leq 
i\leq n-1$}},\\ c_n+2c_{n+1}+\cdots+2c_{2n-1}& {\text{if $i=n$}},  \\
 c_{n-j+1}+\cdots+c_{2n-1}&{\text{if $1\leq i-n=j\leq 
n-1$}},\end{cases}\end{equation*}
if $\frak{g}$ is of type $B_n$;
\begin{equation*} 
k_i =\begin{cases} c_i+c_{2n-i+1}+c_{2n-i+2}+\cdots+c_{2n-1}& {\text{if $2\leq 
i\leq n-1$}},\\ 2c_n+c_{n+1}+\cdots+c_{2n-1}& {\text{if $i=n+1$}},\\
c_{n-j+1}+\cdots+c_{n-1}+2c_n+c_{n+1}+\cdots+ c_{2n-1}&{\text{if $2\leq 
i-n=j\leq n-1$}},
\end{cases}\end{equation*}
if $\frak{g}$ is of type $C_n$;
\begin{equation*} k_i=\begin{cases} c_1 &{\text{if $i=1$}},\\
c_i+c_{2n-i+1}+c_{2n-i+2}+\cdots+c_{2n-2} &{\text{if
$2\leq i\leq n$}},\\
c_{n-j}+\cdots+c_{2n-2} &{\text{if $1\leq i-n=j\leq
n-2$}}.\end{cases} \end{equation*}
if $\frak{g}$ is of type $D_n$.

\medskip
The formulae for the exceptional algebras are naturally much more complicated to
write down. We give as an example the case of $E_6$, and omit the other cases.
If $\frak{g}$ is of type $E_6$, then
\begin{equation*}k_1=c_9,\ \
k_2=c_5+c_{9},\ \
k_3=c_3+c_{5} +c_{9}, \end{equation*}
\begin{equation*} k_4=c_4+c_6+c_8+c_{10}+c_{16},\ \
k_5=c_2+c_3+c_5+c_9, \ \ k_6=c_4+c_7+c_{11}\end{equation*}
\begin{equation*} 
k_7 =c_6+c_7+c_{12}+c_{14}+c_{15}+c_{16},\ \
k_8 =c_{10}+c_{11}+\cdots+c_{15}+c_{16},\ \ \end{equation*}
\begin{equation*}
k_9 =c_1+c_2+c_{3}+c_{5}+c_{9},\ \ 
k_{10}=c_4+c_{8}+c_{13}+c_{14}+c_{15}+c_{16},\ \end{equation*}
\begin{equation*}
k_{11}=c_6+c_{8}+c_{12}+c_{13}+c_{14}+c_{15}+c_{16},\end{equation*}
\begin{equation*}
k_{12}=c_7+c_8+c_{11}+c_{12}+c_{13}+c_{14}+c_{15}+c_{16},\end{equation*}
\begin{equation*}
k_{13}=c_8+c_{10}+c_{11}+c_{12}+c_{13}+c_{14}+c_{15}+c_{16}, \end{equation*}
\begin{equation*}
k_{14}=
c_7+c_8+c_{10}+c_{11}+c_{12}+c_{13}+c_{14}+c_{15}+c_{16},\end{equation*}
\begin{equation*}
k_{15}=c_6+c_{7}+c_{8}+c_{10}+c_{11}+c_{12}+c_{13}+c_{14}+c_{15}+c_{16},
\end{equation*} 
\begin{equation*}
k_{16}=c_4+c_{7}+c_{8}+c_{10}+c_{11}+c_{12}+c_{13}+c_{14}+c_{15}+c_{16}.
\end{equation*}
\medskip

As two examples we write down the details of proof of Lemma 2.3 for type $A_n, 
B_n$. We shall
write $x\equiv y\mod (>\boc)$ if $x-y$ is a $\bold Z[q,q^{-1}]$-linear 
combination of
$E_\bod\ (\bod>\boc)$.

For type  $A_n$  we denote by $E_{i,j}$ ($i>j$) the root vector corresponding to 
$\a_i+\a_{i-1}+\cdots+\a_j$. We have

$$\begin{array}{rl}
M_{f(\boc)}&=E_1^{(c_n)}E_2^{(c_n+c_{n-1})}\cdots
E_n^{(c_n+\cdots+ c_1)}\\[3mm]
&=E_1^{(c_n)}E_2^{(c_n+c_{n-1})}\cdots
E_{n-2}^{(c_n+\cdots+ c_3)}\ds_{\stackrel{\scriptstyle r_n, s_n,t_n \in \bold 
N}{\stackrel{\scriptstyle
r_n+s_n= c_n+\cdots+c_1}    
{s_n+t_n=c_n+\cdots+c_2 } } } v^{- r_n t_n} E^{(r_n)}_{n} E^{(s_n)}_{n,n-1} 
E^{(t_n)}_{n-1}
\\[5mm]
&=\ds_ {\stackrel{\scriptstyle r_n, s_n,t_n \in \bold N }  { 
\stackrel{\scriptstyle r_n+s_n= c_n+\cdots+c_1}    
{s_n+t_n=c_n+\cdots+c_2 } } } v^{- r_n t_n}E_1^{(c_n)}E_2^{(c_n+c_{n-1})}\cdots
E_{n-2}^{(c_n+\cdots+ c_3)} E^{(r_n)}_{n} E^{(s_n)}_{n,n-1} 
E^{(t_n)}_{n-1}\\[5mm]&=\ds _ {\stackrel{\scriptstyle r_n, s_n,t_n \in \bold N}  
  {\stackrel{\scriptstyle r_n+s_n= c_n+\cdots+c_1}    
{s_n+t_n=c_n+\cdots+c_2 } } } v^{- r_n t_n 
}E^{(r_n)}_{n}E_1^{(c_n)}E_2^{(c_n+c_{n-1})}\cdots
E_{n-2}^{(c_n+\cdots+ c_3)}  E^{(s_n)}_{n,n-1} E^{(t_n)}_{n-1}\\[5mm]
&\equiv E_n^{(c_1)}E_1^{(c_n)}E_2^{(c_n+c_{n-1})}\cdots
E_{n-2}^{(c_n+\cdots+ c_3)}  E^{(c_n+\cdots+c_2)}_{n,n-1}   \mod   
(>\boc)\\[5mm]
&=E_1^{(c_n)}E_2^{(c_n+c_{n-1})}\cdots
E_{n-3}^{(c_n+\cdots+ c_4)}\\[5mm]
&\qquad\times \ds _ {\stackrel{\scriptstyle r_{n-1}, s_{n-1},t_{n-1} \in \bold N 
}   {\stackrel{\scriptstyle r_{n-1}+s_{n-1}= c_n+\cdots+c_2}    
{s_{n-1}+t_{n-1}=c_n+\cdots+c_3 } } } v^{- r_{n-1} t_{n-1}} 
E^{(r_{n-1})}_{n,n-1} E^{(s_{n-1})}_{n,n-2} E^{(t_{n-1})}_{n-2}\\[5mm]
&\equiv E_n^{(c_1)}E_{n,n-1}^{(c_2)}E_1^{(c_n)}E_2^{(c_n+c_{n-1})}\cdots
E_{n-3}^{(c_n+\cdots+ c_4)}  E^{(c_n+\cdots+c_3)}_{n,n-2}   \mod   
(>\boc)\\[5mm]
& ...\\[5mm]
&\equiv E_n^{(c_1)}E_{n,n-1}^{(c_2)}\cdots 
E_{n,3}^{(c_{n-2})}E_1^{(c_n)}E_{n,2}^{(c_n+c_{n-1})}   \mod   (>\boc)\\[5mm]
&=E_n^{(c_1)}E_{n,n-1}^{(c_2)}\cdots E_{n,3}^{(c_{n-2})}\ds _ 
{\stackrel{\scriptstyle r_2, s_2,t_2 \in \bold N}    {\stackrel{\scriptstyle 
r_2+s_2= c_n+c_{n-1}}    
{s_2+t_2=c_n } } } v^{- r_2 t_2} E^{(r_2)}_{n,2} E^{(s_n)}_{n,1} 
E^{(t_n)}_{1}\\[5mm]
&\equiv E_n^{(c_1)}E_{n,n-1}^{(c_2)}\cdots 
E_{n,3}^{(c_{n-2}}E_{n,2}^{(c_{n-1})}E_{n,1}^{(c_n)}   \mod   (>\boc).
\end{array}$$

This proves Lemma 2.3 for type $A_n$. 

For type $B_n$, we shall write $E_{i,j}$ $(i\geq j)$ for the root vector 
corresponding to $\a_i+\a_{i-1}+\cdots+\a_j$ and  $E'_{n,j}$ $(n>j)$ for the 
root vector corresponding to
 $\a_n+\a_{n-1}+\cdots+\a_{j+1}+2\a_{j}+\cdots+2\a_1$.
Note that $\a_2,...,\a_n$ generate a root system of type $A_{n-1}$.
 
$$\begin{array}{rl}
M_{f(\boc)}&=E_n^{(k_1)}E_{n-1}^{(k_2)}\cdots E_1^{(k_n)}\cdots
E_n^{(k_{2n-1})}\\[3mm]
&\equiv E_n^{(k_1)}E_{n-1}^{(k_2)}\cdots E_1^{(k_n)}E_2^{(k_{n+1})}\cdots 
E_{n-2}^{(k_{2n-3})}E_n^{(c_2)}E_{n,n-1}^{(k_{2n-2})}  \mod ( >\boc)\\[3mm]

&\equiv E_n^{(k_1)}E_{n-1}^{(k_2)}E_n^{(c_2)}E_{n-2}^{(k_3)}\cdots 
E_1^{(k_n)}E_2^{(k_{n+1})}\cdots E_{n-2}^{(k_{2n-3})}E_{n,n-1}^{(k_{2n-2})}  
\mod ( >\boc)\\[3mm]

&\equiv E_n^{(c_1)}E_{n,n-1}^{(c_2)}E_{n-1}^{(c_{2n-1})}E_{n-2}^{(k_3)}\cdots 
E_1^{(k_n)}E_2^{(k_{n+1})}\cdots E_{n-2}^{(k_{2n-3})}\\[3mm]
&\qquad\times E_{n,n-1}^{(k_{2n-2})}  \mod ( >\boc)\\[3mm]&\equiv 
E_n^{(c_1)}E_{n,n-1}^{(c_2)}E_{n-1}^{(c_{2n-1})}E_{n-2}^{(k_3)}\cdots 
E_1^{(k_n)}E_2^{(k_{n+1})}\cdots E_{n-3}^{(k_{2n-4})}\\[3mm]
&\qquad\times E_{n,n-1}^{(c_3)} E_{n,n-2}^{(k_{2n-3})}  \mod ( >\boc)\\[3mm]

&\equiv 
E_n^{(c_1)}E_{n,n-1}^{(c_2)}E_{n,n-2}^{(c_3)}E_{n-1}^{(c_{2n-1})}E_{n-2}^{(c_{2n
-2}+c_{2n-1})}E_{n-3}^{(k_4)}\cdots E_1^{(k_n)}\\[3mm]
&\qquad\times E_2^{(k_{n+1})}\cdots E_{n-3}^{(k_{2n-4})} E_{n,n-2}^{(k_{2n-3})}  
\mod ( >\boc)\\[3mm]
&...\\
&\equiv E_n^{(c_1)}E_{n,n-1}^{(c_2)}\cdots 
E_{n,2}^{(c_{n-1})}E_{n-1}^{(c_{2n-1})}E_{n-2}^{(c_{2n-2}+c_{2n-1})}\cdots 
E_2^{(c_{n+2}+\cdots+c_{2n-1})}E_1^{(k_n)}\\[3mm]
&\qquad \times E_{n,2}^{(c_{n}+\cdots+c_{2n-1})}  \mod ( >\boc)\\[3mm]
&\equiv E_n^{(c_1)}E_{n,n-1}^{(c_2)}\cdots 
E_{n,2}^{(c_{n-1})}E_{n-1}^{(c_{2n-1})}E_{n-2}^{(c_{2n-2}+c_{2n-1})}\cdots 
E_2^{(c_{n+2}+\cdots+c_{2n-1})}E_{n,1}^{(c_n)}\\[3mm]
&\qquad \times{E'}_{n,1}^{(c_{n+1}+\cdots+c_{2n-1})}  \mod ( >\boc)\\[3mm]
&\equiv E_n^{(c_1)}E_{n,n-1}^{(c_2)}\cdots 
E_{n,2}^{(c_{n-1})}E_{n,1}^{(c_n)}E_{n-1}^{(c_{2n-1})}E_{n-2}^{(c_{2n-2}+c_{2n-1
})}\cdots E_2^{(c_{n+2}+\cdots+c_{2n-1})}\\[3mm]
&\qquad\times {E'}_{n,1}^{(c_{n+1}+\cdots+c_{2n-1})}  \mod ( >\boc)\\[3mm]
&\equiv E_n^{(c_1)}E_{n,n-1}^{(c_2)}\cdots 
E_{n,2}^{(c_{n-1})}E_{n,1}^{(c_n)}E_{n-1}^{(c_{2n-1})}E_{n-2}^{(c_{2n-2}+c_{2n-1
})}\cdots E_3^{(c_{n+3}+\cdots+c_{2n-1})}\\[3mm]
&\qquad\times {E'}_{n,1}^{(c_{n+1})}{E'}_{n,2}^{(c_{n+2}+\cdots+c_{2n-1})}  \mod 
( >\boc)\\[3mm]
&\equiv E_n^{(c_1)}E_{n,n-1}^{(c_2)}\cdots 
E_{n,2}^{(c_{n-1})}E_{n,1}^{(c_n)}{E'}_{n,1}^{(c_{n+1})}E_{n-1}^{(c_{2n-1})}E_{n
-2}^{(c_{2n-2}+c_{2n-1})}\cdots \\[3mm]
&\qquad\times E_3^{(c_{n+3}+\cdots+c_{2n-1})} 
{E'}_{n,2}^{(c_{n+2}+\cdots+c_{2n-1})}  \mod ( >\boc)\\[3mm]
&...\\[3mm]
&\equiv E_n^{(c_1)}E_{n,n-1}^{(c_2)}\cdots 
E_{n,2}^{(c_{n-1})}E_{n,1}^{(c_n)}{E'}_{n,1}^{(c_{n+1})}\cdots 
{E'}_{n,n-1}^{(c_{2n-1})}  \mod ( >\boc)
.\end{array}$$

For other types the proof is completely similar.

Finally we conjecture for each reduced expression of $w_0$ there
is a corresponding monomial basis of $\bu^+$.

\bigskip
\noindent{\bf Acknowledgement:} N.X.  was  supported  by the
National Natural Science Foundation of China. N.X. also would like
to thank Newton Institute  at Cambridge University for financial
support and hospitality during his visit to the institute.

\end{document}